\documentclass[11pt]{article} \usepackage[a4paper,scale=0.7]{geometry} \usepackage{amsmath,amssymb,amsthm,mathtools} 
\usepackage{graphicx} 
\usepackage[T1]{fontenc} 
\usepackage{lmodern}
\usepackage{hyperref} 
\newtheorem{theorem}{Theorem} 
\newtheorem{proposition}[theorem]{Proposition} 
 
\newtheorem{lemma}[theorem]{Lemma} 
\theoremstyle{remark} 
\newtheorem{remark}[theorem]{Remark}
\AtBeginEnvironment{remark}{%
  \pushQED{\qed}%
}
\AtEndEnvironment{remark}{\popQED}
\usepackage{dsfont}

\usepackage{todonotes} 
 
\usepackage{libertine}

\usepackage[maxnames=99, backend=bibtex, url=false, isbn=false, doi=false, style=alphabetic]{biblatex} 
\renewbibmacro{in:}{} 
\newcommand{\T}{\mathds T} 
\newcommand{\E}{\mathbb E} 
\newcommand{\Pp}{\mathbb P} 
\newcommand{\1}{\mathbf 1} 
\newcommand{\bb}{\mathbf b} 
\newcommand{\bn}{\mathbf n} 

\usepackage[affil-sl]{authblk}
\title{The local survival threshold of the contact process\\ on regular and periodic trees} 

\usepackage{orcidlink}
\author[ ]{John Fernley \orcidlink{0000-0002-6635-4341}}
\affil[ ]{University of Warwick, Department of Statistics, CRiSM}
\affil[ ]{\rm john.fernley@warwick.ac.uk}
\date{\today}

\date{\today}

\begin{document}
\maketitle

\begin{abstract}
We consider the contact process on quasi-transitive trees whose degrees repeat with a fixed period along each line of descent, including the regular tree. As the degrees grow in asymptotically fixed proportions, we prove that the local-survival critical value $\lambda_2$, multiplied by the square root of their geometric mean, converges to a limit $\kappa$. In the limiting root process, Poisson edges on the root's timeline represent possible returns of infection through the branches. The rescaled limit $\kappa$ is the unique rate $a>0$ for which the geometric mean of the expected total infected times, one for each root type, equals $1/a$ in the limit model. We also give strict upper and lower bounds on $\kappa$ for every period. On the regular tree our characterisation answers a question of Liggett, and we show that $\kappa<2/3$, improving his upper bound of $1$.
\par\smallskip
\noindent\textit{2020 Mathematics Subject Classification:} Primary 60K35, secondary 05C05, 82C22.
\end{abstract}

\section{Local survival of the contact process}\label{sec:contact_intro}

The contact process on a graph $G=(V,E)$ with infection rate $\lambda>0$ is a continuous-time Markov process $(\xi_t)_{t\ge 0}$ with values in $\{0,1\}^V$ describing the spread of an infection on $G$, with the interpretation that $\xi_t(x)=1$ if the individual $x$ in $V$ is infected at time $t$, while $\xi_t(x)=0$ if it is healthy. We identify the configuration with its set of infected vertices and write $\xi_t^A$ for the process started with precisely the vertices of $A$ infected. The dynamics of the process are that each infected vertex recovers at rate 1, and each edge between an infected vertex and a healthy one transmits the infection to the healthy vertex at rate $\lambda$.

The contact process on a regular or periodic tree can survive globally while visiting each fixed vertex only finitely often. Starting from a single infected root $\circ\in V$, we say that the process survives locally if the root is infected at arbitrarily late times. Write $\lambda_2$ for the critical infection rate for local survival.

How can the infection return to the root after it has recovered? On a $d$-regular tree at infection rate $\lambda=a/\sqrt d$, the probability that a given neighbour offers a route out and back during a fixed time interval is of order $1/d$. There are $d$ neighbours, and in the limit these possible returns form a Poisson process of pairs of times $(s,t)$. Such a pair can carry infection from the root at time $s$ through a single branch and back to the root at time $t$, even if the root has recovered in between. We show that the root's expected occupation time in this limiting process determines the limit of $\sqrt d\,\lambda_2$.

For the contact process on the tree started from the root alone, the probabilities of global survival and local (or strong) survival are respectively
\begin{align*}
s_1(\lambda)&=\Pp_\lambda\big(\, \forall t\geq 0, \; \xi_t^{\{\circ\}}\neq \varnothing
\big),\\
s_2(\lambda)&=\Pp_\lambda\big(\, \xi_t^{\{\circ\}}(\circ)=1 \text{ for arbitrarily large }t \big)\le s_1(\lambda).
\end{align*}

Because $s_1$ and $s_2$ are necessarily non-decreasing in $\lambda$ we can define the critical parameters
\begin{align*}
\lambda_1&=\sup\left\{ \lambda : s_1(\lambda)=0
\right\},\\
\lambda_2&=\sup\left\{ \lambda : s_2(\lambda)=0 \right\},
\end{align*}
and we automatically have $\lambda_1\le \lambda_2$. In \cite{MR4125782} it was found that on periodic trees in fact $0<\lambda_1< \lambda_2<\infty$, generalising the same statement on $d$-regular trees which is due to \cite{MR1188054} for $d\ge4$ and \cite{MR1415225} for $d=3$ (note for the line $d=2$ we have instead  $0<\lambda_1= \lambda_2<2$). Neither critical value depends on the chosen root. Indeed, infection can be transmitted along any prescribed finite path with positive probability; the Markov property and repeated such trials transfer local survival between vertices.

For regular trees of degree $d$, it is known that $\lambda_1\sim 1/d$. Liggett \cite[after Equation 5.5]{MR1717346} asked for the value of the rescaled limit $\kappa(1)=\lim_{d\to\infty}\sqrt d\,\lambda_2$. For fixed $a>0$, write $q(t;a)$ for the probability that the initially infected root is infected at time $t$ in the limiting process. The Poisson pairs have intensity $a^2q(t-s;a)\,ds\,dt$ for $0<s<t$, with $a^2=d\lambda^2$. A return through a neighbour requires an infection path within its branch from that neighbour at time $s$ to the same vertex at time $t$. Since the process at the branch root has the same limiting law as at the original root, the limiting probability of such a path is $q(t-s;a)$. Together with independent rate-one recoveries, this rule determines $q$ uniquely for each $a$ (Lemma~\ref{lem:fixed-point}). By Theorem~\ref{thm:periodic}, $\kappa(1)$ is the unique $a>0$ satisfying
\begin{equation}\label{eq:regular-constant}
a\int_0^\infty q(t;a)\,dt=1.
\end{equation}

Pemantle~\cite{MR1188054} gave asymptotic lower and upper bounds of $2-\sqrt2$ and $e$ for $\sqrt d\,\lambda_2$. The lower bound follows from
\[
\lambda_2 \geq \frac{2}{2\sqrt{d-1}-1+\sqrt{2d-3}}.
\] 

Liggett~\cite{MR1717346} improved the asymptotic upper bound from $e$ to $1$ by showing that
\[
\lambda_2 \leq \frac{1}{\sqrt{d-1}-1},
\]
and in Proposition \ref{prop:upper-bounds} we show here that in fact $\kappa(1)<2/3$ (and the same bound applies to all quasi-transitive periodic trees). This understanding of $\lambda_2$ is translated moreover to the random regular graph in the appearance of a logarithmic length lapse in the surviving epidemic \cite{MR5091097}.

Later in this article we will generalise \eqref{eq:regular-constant} to the periodic trees of the following section, using one expected root occupation time for each type in the limiting model, by which we can get a general sharp asymptotic for offspring numbers growing in asymptotic positive proportions. Highly unequal degree regimes, in which one offspring number is polynomially larger than the others, are treated by \cite{MR4089731,HuangDurrett2023}, where the second critical value is found on a different order, so local survival in generality for periodic trees remains a complicated issue.

\section{The large-degree periodic tree}

An \emph{end} of a tree is an equivalence class of infinite simple paths (rays), where two rays are equivalent if they eventually coincide. For fixed $k\ge1$ and positive integers $\bn=(n_1,\ldots,n_k)$ with $\prod_i n_i>1$, let $\T(\bn)$ be the periodic tree with one \emph{distinguished end} chosen as part of its structure. From each vertex the first edge on the unique ray to this end leads to its parent. A type-$i$ vertex has a parent of type $i-1$ and $n_i$ children of type $i+1$, its neighbours in the other direction. This is the quasi-transitive version of the periodic tree \cite{MR1436761,MR4125782}. There is no exceptional boundary vertex as in \cite{MR1062053,MR4089731}, every vertex has a parent, and degree $d_i=n_i+1$ at type $i$.

Two vertices lie in the same orbit if a graph automorphism carries one onto the other. A graph is \emph{transitive} if it has one vertex orbit and \emph{quasi-transitive} if it has finitely many. In $\T(\bn)$, any two vertices of the same type can be mapped to one another by matching their rays towards the distinguished end and the branches attached along those rays. Thus $\T(\bn)$ has at most $k$ vertex orbits. A regular tree has one orbit, while a biregular tree with two distinct degrees has exactly two, one for each degree.

Consider a sequence $\bn^{(N)}=(n_{1,N},\ldots,n_{k,N})$ of degree vectors of length $k$. Set the product and normalised version
\begin{equation}\label{eq:normalization}
\Gamma_N=\prod_{i=1}^k n_{i,N},\qquad b_{i,N}=\frac{n_{i,N}}{\Gamma_N^{1/k}}.
\end{equation}

Throughout the limit arguments we assume
\begin{equation}\label{eq:shape-assumption}
\Gamma_N\longrightarrow\infty,\qquad b_{i,N}\longrightarrow b_i>0\quad(1\le i\le k).
\end{equation}

Note $\prod_i b_i=1$ so the vector is normalised with respect to the geometric mean, and every offspring number tends to infinity. We call $\bb=(b_1,\ldots,b_k)$ the normalised degree vector. The period $k$ is fixed.

\section{Results}

The limiting constant is determined by a family of root processes, one for each vertex type. For $c>0$ and measurable $f:[0,\infty)\to[0,1]$, let $X^{c,f}(t)\in\{0,1\}$ be a root process started from state $1$ (infected). Place rate-one Poisson recovery marks on the timeline and independent Poisson edges $(s,w)$ with $0<s<w$ and intensity
\begin{equation}\label{eq:edge-intensity}
c^2 f(w-s)\,ds\,dw.
\end{equation}

A recovery sets the state to $0$. An edge $(s,w)$ sets the state at $w$ to $1$ if the root was infected at $s$ and otherwise leaves it unchanged. Write $(\Psi_cf)(t)=\Pp(X^{c,f}(t)=1)$.

At rescaled infection rate $a>0$, the process at a type-$i$ root is $X^{a\sqrt{b_i},q_{i+1}}$, where $q_{i+1}$ is the infection probability profile for the next type. The profiles therefore solve
\begin{equation}\label{eq:cyclic-root}
q_i=\Psi_{a\sqrt{b_i}}q_{i+1},\qquad 1\le i\le k,
\end{equation}
with $q_{k+1}:=q_1$. Lemma~\ref{lem:fixed-point} proves existence and uniqueness.

We call the expected total infected time of an initially infected type-$i$ root in the limiting process its \emph{root susceptibility}. Integrating its infection probability $q_i(t)$ gives
\begin{equation}\label{eq:H}
\chi_i(a,\bb)=\int_0^\infty q_i(t)\,dt,\qquad H_{\bb}(a)=a^k\prod_{i=1}^k\chi_i(a,\bb).
\end{equation}

\begin{theorem}\label{thm:periodic}
For every $k\geq 1$ and normalised vector $\bb\in\mathds{R}_+^k$ there is a constant $\kappa(\bb)\in(0,\infty)$ such that, under \eqref{eq:shape-assumption},
\begin{equation}\label{eq:periodic-limit}
\Gamma_N^{1/(2k)}\lambda_2(\T(\bn^{(N)})) \longrightarrow\kappa(\bb).
\end{equation}

The constant $\kappa(\bb)$ is the unique $a>0$ satisfying
\begin{equation}\label{eq:kappa}
H_{\bb}(a)=1.
\end{equation}
\end{theorem}

Pemantle \cite[after Equation 7]{MR1188054} gave the upper bound $\lambda_2=O(\Gamma_N^{-1/(2k)})$ for rooted periodic trees, which applies automatically to our quasi-transitive version in which $\circ$ has an extra neighbour. To answer some ambiguous conjectures in that work, \cite{HuangDurrett2023} recently gave sharp asymptotics for $\lambda_2$ on a smaller order than $\Gamma_N^{-1/(2k)}$, with the assumption that degrees in the periodic tree are unbalanced: one degree goes to infinity and all others are polynomially sublinear by comparison.

Our asymptotic instead makes the fixed-proportion assumption \eqref{eq:shape-assumption}, so the result is on the same order as Pemantle's upper bound, and this also allows us to include the regular tree in the theorem. For biregular trees $\T(n,m)$ with $m=n^{1-\Theta(1)}$, Huang and Durrett's asymptotic expression
\[
\lambda_2(\T(n,m))\sim \sqrt{ \frac{\log(n/m)}{2n} }
\]
also arises when we first let the degrees grow with a fixed ratio and then let their limiting ratio tend to infinity (Remark~\ref{rem:biregular-comparison}).

We also present bounds on $\kappa(\bb)$ which use a symmetric $k\times k$ matrix $M_{\bb}$. For $k\ge3$, it has $\sqrt{b_1},\ldots,\sqrt{b_{k-1}}$ immediately above and below the diagonal. The weight $\sqrt{b_k}$ joins the last type to the first and appears in the upper right and lower left corners. With blank entries denoting zeros, the matrix is
\begin{equation}\label{eq:cyclic-matrix}
\renewcommand{\arraystretch}{1.2} \setlength{\arraycolsep}{2pt} \newcommand{\matrixcell}[1]{\mathmakebox[2.8em]{#1}} M_{\bb}=
 \begin{pmatrix}
\matrixcell{0} & \matrixcell{\sqrt{b_1}} & \matrixcell{} &
  \matrixcell{} & \matrixcell{} & \matrixcell{\sqrt{b_k}} \\
  \sqrt{b_1} & 0          & \sqrt{b_2} &        &                &                \\
             & \sqrt{b_2} & 0          & \ddots &                &                \\
             &            & \ddots     & \ddots & \ddots         &                \\
             &            &            & \ddots & 0              & \sqrt{b_{k-1}} \\
\sqrt{b_k} &            &            &        & \sqrt{b_{k-1}} & 0
 \end{pmatrix}.
\end{equation}

Section~\ref{sec:brw-threshold} explains the matrix weights in terms of the branching random walk. For $k=1,2$ the preceding and following types coincide, so their weights add, giving expressions
\[
M_{(1)}=(2),\qquad M_{(b_1,b_2)}=
 \begin{pmatrix}
  0                     & \sqrt{b_1}+\sqrt{b_2}\\
\sqrt{b_1}+\sqrt{b_2} & 0
 \end{pmatrix}.
\]

We write $\rho(M_{\bb})$ for the Perron eigenvalue of this matrix, equivalently the spectral radius.

\begin{proposition}\label{prop:upper-bounds}
For every $k\geq 1$ and normalised vector $\bb\in\mathds{R}_+^k$,
\[
\frac{4-2\sqrt2}{\rho(M_{\bb})} \ <\ \kappa(\bb)\ <\ \frac23.
\]
\end{proposition}

Since $\rho(M_{(1)})=2$, the bounds for the regular tree are
\[
2-\sqrt2< \kappa(1) <\frac{2}{3}.
\]

\begin{figure}[b!]
\centering \includegraphics[width=0.75\textwidth]{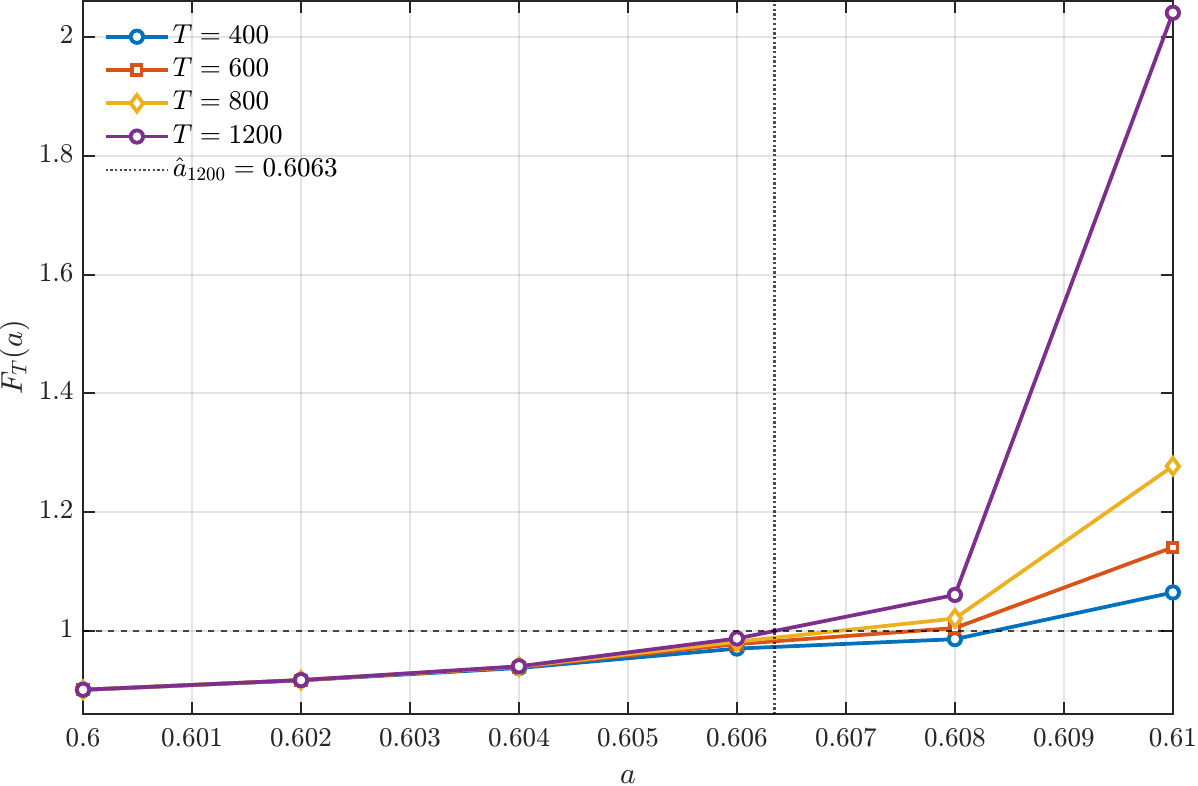} \caption{The simulation suggests that the regular-tree constant \(\kappa(1)=\lim_{d\to\infty}\sqrt d\,\lambda_2\) is about \(0.606\). The curves estimate \(F_T(a)=a\int_0^T q(t;a)\,dt\), and note that the exact solution of \(F_T(a)=1\) approaches \(\kappa(1)\) from above as \(T\to\infty\). At each marked rate \(a\), we started from the profile \(e^{-t}\) and estimated 60 successive iterates of \(\Psi_a\), using \(10^5\) paths per iterate and a grid of spacing \(0.1\) for the delay profile. Linear interpolation between rates spaced by \(0.002\) gives estimated crossings of the dashed line \(F_T(a)=1\) at \(0.6083,\,0.6076,\,0.6069,\,0.6063\) for \(T=400,\,600,\,800,\,1200\), respectively.}\label{figure_sim}
\end{figure}

The lower end agrees with Pemantle's asymptotic bound \cite[Theorem 2.2]{MR1188054}. Reducing Liggett's upper bound \cite{MR1717346} from $1$ to $2/3$ leaves an interval less than one fifth as wide. Simulation of the limiting root process suggests $\kappa(1)\approx 0.606$ (Figure~\ref{figure_sim}).

\begin{remark}[Branching random walk comparisons]\label{rem:brw-comparison}
For the branching random walk, the rescaled local survival threshold converges to $1/\rho(M_{\bb})$, as shown in Section~\ref{sec:brw-threshold}. Domination by the branching random walk gives the comparison
\begin{equation}\label{eq:main-bounds}
\frac1{\rho(M_{\bb})}\le\kappa(\bb).
\end{equation}

For period $2$ this is Pemantle's bound \cite[before Equation~8]{MR1188054}, but of course Proposition~\ref{prop:upper-bounds} improves this lower bound for every $\bb$. For the limiting ratio $\theta$ of $n_{2,N}/n_{1,N}$, the branching random walk constant $(\theta^{-1/4}+\theta^{1/4})^{-1}$ is asymptotic to $\theta^{-1/4}$ as $\theta\to\infty$. At the point of local survival, the exclusion makes a significant difference and so the contact process constant is larger by a factor asymptotic to $\sqrt{(\log\theta)/2}$, as in the following remark. This is different to the global survival threshold, where the exclusion is not significant and the critical values are asymptotically the same \eqref{eq:global-asymptotic}.
\end{remark}

\begin{remark}[Biregular trees]\label{rem:biregular-comparison}
Let $n>m$ be the two offspring numbers and write $\lambda_2=\lambda_2(\T(n,m))$. The biregular tree was considered in \cite[Theorem~3.3]{MR1188054}, with $m=1$, where the result is $\lambda_2=O(\sqrt{\log n/n})$ as $n\to\infty$. Our result and \cite{MR4089731} give, respectively,
\begin{align*}
\lim_{\theta\to\infty}\ \lim_{n\to\infty} \frac{2n\lambda_2^2}{\log(n/m)} &=1 &&\text{on }\T(n,n/\theta+o(n)),\\[6pt] \lim_{n\to\infty} \frac{2n\lambda_2^2}{\log(n/m)} &=1 &&\text{on }\T(n,n^{\beta+o(1)}) \text{ with } 0\le\beta<1.
\end{align*}

For the first line, label the type with $m$ offspring as type $1$. At fixed $\theta>1$, Theorem~\ref{thm:periodic} gives the inner limit
\[
\frac{2\sqrt\theta}{\log\theta}\,\kappa(\theta^{-1/2},\theta^{1/2})^2.
\]

Proposition~\ref{prop:biregular-ratio} in Appendix~\ref{app:biregular-ratio} then gives the outer limit. The second line is not exactly a result of \cite{MR4089731} because their trees have a special root vertex with no parent, but it follows easily from sandwiching the quasi-transitive biregular tree between $(n,m)$ and $(n+1,m)$ versions of  their tree.

Because the two limits above are quite different, the proofs are as well: for Huang and Durrett the process is understood by communication between stars of degree on a polynomially larger order, while in this article the degrees are comparable and returns through children of the opposite type give coupled fixed-point equations for the two root infection probabilities.
\end{remark}

\subsection{Outline of the proof}

Set $\lambda=a\Gamma_N^{-1/(2k)}$. On a bounded time interval, the leading contribution from a child comes from a single arrow to the child at time $s$ and an arrow back at time $t$. The pair can carry infection only if the child, of type $i+1$, has an infection path from $s$ to $t$ within its branch. This gives the factor $q_{i+1}(t-s)$ in the Poisson intensity, while $n_{i,N}\lambda^2\to a^2b_i$ gives its prefactor. Cutting the edges of a fixed simple path then expresses its leading infection probability as a convolution of the root profiles.

For extinction, we sum these path contributions with weights adapted to the offspring numbers. A simple path first moves towards the distinguished end and then away from it. Each complete period in either segment contributes a factor of $H_{\bb}(a)$ to the time-integrated limiting weighted mass. When $H_{\bb}(a)<1$, the sum is finite, so there is a common time $T$ at which the limiting expected weighted mass is less than $1$ for every starting type, after normalisation by the starting vertex's weight. By convergence at time $T$ and additivity, for all sufficiently large $N$ the expected weighted mass after time $T$ is at most a fixed fraction of its initial value for every finite initial infection set. Iterating the bound and controlling the mass between time steps gives finite expected infected time at the root and hence local extinction.

For survival, remove the edge from a type-$i$ root to a prescribed child. For fixed $T$, the root's expected occupation time up to $T$ in its component converges to $\int_0^T q_i(t)\,dt$ as $N\to\infty$. Counting arrows to the child while the root is infected gives a lower bound on the probability of hitting that child in the full tree. A branching-and-return argument gives local survival when the product of these $k$ hitting probabilities exceeds $\Gamma_N^{-1/2}$. Thus when $H_{\bb}(a)>1$, the contact process survives locally on the trees for all sufficiently large $N$.

\section{Branching random walk bounds}
\label{sec:prelim}

We use independent Poisson recovery marks of rate $1$ at vertices and infection arrows of rate $\lambda$ on directed edges. Infection follows paths moving forward in time, avoiding recovery marks and using arrows in their direction.

Deleting vertices or arrows can only reduce infection. The construction also gives additivity: the process from $A$ is the union of the processes from its individual vertices using the same marks. On any fixed undirected subgraph $F$, reversal of the graphical construction in a fixed time interval gives self-duality:
\begin{equation}\label{eq:duality}
\Pp_F^{\{x\}}(y\in\xi_t) =\Pp_F^{\{y\}}(x\in\xi_t).
\end{equation}

We will use the following implication repeatedly:
\begin{equation}\label{eq:occupation-extinction}
\int_0^\infty\Pp^A(\circ\in\xi_t)\,dt<\infty \quad\Longrightarrow\quad \text{no local survival at }\circ.
\end{equation}

Indeed, the expected number of recovery marks occurring while $\circ$ is infected equals this integral. A finite expectation gives only finitely many such recoveries almost surely. Each infected period at $\circ$ lasts from infection until the next recovery mark, which occurs in finite time almost surely. Hence there are only finitely many such periods at $\circ$.

We state here a couple of very simple supermartingale-based branching random walk bounds that will be used later.

\begin{lemma}\label{lem:weighted-brw}
Let $w_x>0$ be weights on a bounded-degree graph, with $w_{\circ}=1$, such that
\[
\sum_{y:y\sim x}w_y\le Lw_x\quad\text{for every }x.
\]

Then the contact process started from $\{\circ\}$ satisfies
\begin{align}
\E^{\{\circ\}}\sum_x w_x\1\{x\in\xi_t\}
 &\le e^{(\lambda L-1)t}, \label{eq:brw-total}\\
\E^{\{\circ\}}\sum_{\operatorname{dist}(\circ,x)\ge K} w_x\1\{x\in\xi_t\} &\le e^{-t}\sum_{j\ge K}\frac{(\lambda Lt)^j}{j!}.
 \label{eq:brw-tail}
\end{align}
\end{lemma}

\begin{proof}
Dominate the contact process by a branching random walk: particles die at rate $1$ and give birth along each adjacent edge at rate $\lambda$, allowing multiple particles per vertex. Fix a walk $x_0=\circ,x_1,\ldots,x_j$. Let $m_h(t)$ be the expected number of particles alive at $t$ whose ancestral walk is exactly $x_0,\ldots,x_h$. For $h=0$ only the original particle counts, so $m_0(t)=e^{-t}$. Each particle with ancestral walk $x_0,\ldots,x_{h-1}$ gives birth at $x_h$ at rate $\lambda$, and a child born at $s$ survives to $t$ with probability $e^{-(t-s)}$. Integrating over its birth time and using induction on $h$ gives
\begin{align*}
m_h(t)
 &=\int_0^t \lambda m_{h-1}(s)e^{-(t-s)}\,ds\\
&=\int_0^t \lambda e^{-s}\frac{(\lambda s)^{h-1}}{(h-1)!}
                  e^{-(t-s)}\,ds\\
&=\frac{e^{-t}\lambda^h}{(h-1)!}\int_0^t s^{h-1}\,ds =e^{-t}\frac{(\lambda t)^h}{h!}.
\end{align*}

Summing the path endpoint weights over all $j$-step walks from $\circ$ gives at most $L^j$, by iteration of the assumed row inequality. Summing over $j\ge0$ proves \eqref{eq:brw-total}. A path endpoint at distance at least $K$ requires at least $K$ genealogical steps, giving \eqref{eq:brw-tail}. 
\end{proof}

\subsection{Comparisons to the contact process}\label{sec:brw-threshold}

In the branching random walk, each particle dies at rate $1$ and gives birth at each neighbouring vertex at rate $\lambda$.

For local survival we need to count descendants born back at the starting vertex $\circ$. The initial particle has on average $\lambda$ offspring at each child of $\circ$, and each of those has on average $\lambda$ offspring back at $\circ$. At $\lambda=a\Gamma_N^{-1/(2k)}$ and when $\circ$ has type $i$, the expected number of grandchildren born back at $\circ$ through its $n_{i,N}$ children is therefore
\[
n_{i,N}\lambda^2=a^2 b_{i,N}\longrightarrow (a\sqrt{b_i})^2.
\]

More generally, when $A_N$ denotes the adjacency matrix of the tree, $(A_N^\ell)_{\circ,\circ}$ counts the $\ell$-step walks from $\circ$ to $\circ$. Each step contributes a factor $\lambda$ to the expected number of descendants, so
\[
\E_{\circ}[\text{number of generation $\ell$ particles born at $\circ$}] = \lambda^\ell (A_N^\ell)_{\circ,\circ}.
\]

The tree is bipartite, so returns occur at even lengths. By the return-walk criterion \cite[Lemma~3.1]{MR1880232}, the local critical value of the branching random walk satisfies
\[
\frac{1}{\lambda_2^{\mathrm{BRW}}(\T(\bn^{(N)}))} = \lim_{m\to\infty}\big((A_N^{2m})_{\circ,\circ}\big)^{1/(2m)} = \rho(A_N),
\]
where $\rho(A_N)$ is the tree's adjacency spectral radius.

The grandchildren calculation suggests a weight of $\sqrt{n_{i,N}}$ in each direction between a type-$i$ vertex and a child. For our periodic tree, the theorem of Saloff-Coste and Woess \cite[Theorem~4.3]{SaloffCosteWoess1997} says that the cyclic matrix with these weights has largest eigenvalue $\rho(A_N)$. Writing $\bb_N=(b_{1,N},\ldots,b_{k,N})$, we have $\sqrt{n_{i,N}}=\Gamma_N^{1/(2k)}\sqrt{b_{i,N}}$, so
\[
\frac{\rho(A_N)}{\Gamma_N^{1/(2k)}}=\rho(M_{\bb_N}) \longrightarrow \rho(M_{\bb}).
\]

Suppressing births at occupied vertices gives the contact process, so $\lambda_2^{\mathrm{BRW}}(\T(\bn^{(N)}))\le\lambda_2(\T(\bn^{(N)}))$. Combining the spectral-radius limit with Theorem~\ref{thm:periodic} gives \eqref{eq:main-bounds}.

The gap between the local thresholds reflects the effect of suppressing births at occupied vertices. For global survival the two processes have the same leading-order threshold under \eqref{eq:shape-assumption}. For the contact process,
\begin{equation}\label{eq:global-asymptotic}
\lambda_1(\T(\bn^{(N)})) \sim\Gamma_N^{-1/k}.
\end{equation}

For periods two and three, this was in \cite{JiangEtAl2019}. In general it can be proved by sandwiching the contact process between the process keeping only first infections from parent to child (the SIR infection) and the branching random walk.

\section{The limiting root processes}\label{sec:root}

We first establish existence and uniqueness for the cyclic system
\eqref{eq:cyclic-root}.

Adding edges while keeping the recovery marks fixed shows that $(\Psi_cf)(t)$ is nondecreasing in both $c$ and $f$. If there is no recovery by $t$, the root is still infected, so $(\Psi_cf)(t)\ge e^{-t}$.

\begin{lemma}\label{lem:fixed-point}
For every $a>0$ and positive $\bb$, \eqref{eq:cyclic-root} has a unique solution with values in $[0,1]^k$. Its profiles satisfy
\begin{equation}\label{eq:profile-basic}
e^{-t}\le q_i(t)\le1,
\end{equation}
and each profile is nondecreasing in every coefficient $a^2b_j$. For fixed $\bb$ and $t$, each $q_i(t)$ is continuous in $a$.
\end{lemma}

\begin{proof}
Use the same recovery marks for $\Psi_c f$ and $\Psi_d h$, and couple their common edges. The root states agree at time $t$ unless an excess edge ends by then. The expected number of such edges gives the bound
\begin{equation}\label{eq:volterra}
|(\Psi_c f)(t)-(\Psi_d h)(t)| \le \int_0^t(t-u)|c^2 f(u)-d^2 h(u)|\,du.
\end{equation}

Start with $f_i^{(0)}(t)=e^{-t}$ and $g_i^{(0)}(t)=1$, and iterate
\[
f_i^{(j+1)}=\Psi_{a\sqrt{b_i}}f_{i+1}^{(j)},\qquad g_i^{(j+1)}=\Psi_{a\sqrt{b_i}}g_{i+1}^{(j)}.
\]

Every output at time $t$ lies between $e^{-t}$ and $1$, so by monotonicity the lower iterates increase and the upper iterates decrease.

Fix $A>0$ and put $C=A^2\max_i b_i$. For $0<a\le A$ we claim that
\[
0\le g_i^{(j)}(t)-f_i^{(j)}(t)\le\frac{C^j t^{2j}}{(2j)!}.
\]

The initial gap is at most $1$. If the bound holds at step $j$, then \eqref{eq:volterra} gives
\[
g_i^{(j+1)}(t)-f_i^{(j+1)}(t) \le\frac{C^{j+1}}{(2j)!}\int_0^t(t-u)u^{2j}\,du =\frac{C^{j+1}t^{2j+2}}{(2j+2)!}.
\]

For $t\le T$ the gap is at most $C^jT^{2j}/(2j)!\to0$. Thus both iterations converge to the same profiles $q_i$, uniformly for $t\le T$ and $0<a\le A$. By \eqref{eq:volterra} we may pass to the limit in the recursion for $f_i^{(j+1)}$, so $q$ solves \eqref{eq:cyclic-root}. Any solution lies between the initial bounds and hence, by monotonicity, between the two iterates at every step. It therefore equals $q$.

Increasing any coefficient $a^2b_j$ can only increase the lower iterates and hence their limits. For fixed $\bb$ and $t$, bounded convergence in \eqref{eq:volterra} gives continuity of $f_i^{(j)}(t)$ in $a$ by induction
on $j$. The convergence to $q_i(t)$ is uniform for $0<a\le A$, so $q_i(t)$ is continuous in $a$ as well.
\end{proof}

Each $\chi_i$ is at least $1$ (and possibly $\infty$) and is nondecreasing in $a$, by Lemma~\ref{lem:fixed-point}. This is a susceptibility in the usual sense, but only on the limiting binary-state model, not to be confused with the usual susceptibility of the contact process on a graph which sums the expected occupation times over all vertices, e.g. \cite[Eq.~(1.4)]{AizenmanJung2007}. The extinction and survival arguments below show that $H_{\bb}(a)<1$ and $H_{\bb}(a)>1$ give the two sides of the limiting critical value.

On each bounded time interval, the root infection probabilities for the contact process on the trees converge uniformly to those of the limiting model.

\begin{lemma}\label{lem:root-limit}
Fix $a>0$ and set ${\lambda}=a\Gamma_N^{-1/(2k)}$. For a vertex $\circ$ of type $i$, let $q_{i,N}(t)$ be the probability that $\circ$ is infected at time $t$ when only $\circ$ is initially infected on the full tree $\T(\bn^{(N)})$. Let $(q_i)_{i=1}^k$ be the solution of \eqref{eq:cyclic-root}. For every $T>0$,
\begin{equation}\label{eq:root-limit}
\max_{1\le i\le k}\sup_{0\le t\le T} |q_{i,N}(t)-q_i(t)|\longrightarrow0.
\end{equation}

The same limit holds if any fixed number of edges incident to the initially infected vertex are deleted and the process is run in its remaining component.
\end{lemma}

\begin{proof}
For any vertex $v$, let $S(v)$ consist of $v$ and its descendants. If $v$ has type $j$, write
\[\bar q_{j,N}(u)=\Pp_{S(v)}^{\{v\}}(v\in\xi_u).\]

We first follow the infection at $\circ$ in $S(\circ)$, starting with only $\circ$ infected. The edge from $\circ$ to its parent is absent throughout this calculation. Each of the $n_{i,N}$ children $x$ of $\circ$ has type $i+1$ and is attached to $\circ$ by the single edge $\{\circ,x\}$.

For each child $x$, record a pair $(s,t)$ if the only arrows on $\{\circ,x\}$ in $[0,T]$ are $\circ\to x$ at $s$ and $x\to \circ$ at $t>s$, and there is an infection path inside $S(x)$ from $(x,s)$ to $(x,t)$. To test this last condition, start with $x$ alone infected at $s$ in $S(x)$ and check whether $x$ is infected at $t$. Given the two arrows, the test succeeds with probability $\bar q_{i+1,N}(t-s)$. Thus we can record the pairs independently for different children, before examining the recovery marks at $\circ$.

In the process started from $\circ$, every $S(x)$ starts empty. If $\{\circ,x\}$ has at most two arrows, it can return infection to $\circ$ only through a recorded pair $(s,t)$ for which $\circ$ was infected at $s$. Consequently the recorded pairs, used as edges on the timeline of $\circ$ with its original recovery marks, reproduce the infection history at $\circ$ unless some $\{\circ,x\}$ has at least three arrows. Each $\{\circ,x\}$ has a Poisson number of arrows with mean $2\lambda T$, counting both directions, so
\begin{equation}\label{eq:three-arrows}
\Pp(\text{some $\{\circ,x\}$ has at least three arrows in $[0,T]$}) \le \frac{n_{i,N}(2\lambda T)^3}{6}.
\end{equation}

For one child the probability of a recorded pair in $ds\,dt$ is
\begin{equation}\label{eq:eligible-density}
\lambda^2e^{-2\lambda T}\bar q_{i+1,N}(t-s)\,ds\,dt, \qquad 0<s<t<T.
\end{equation}

The factor $\lambda^2e^{-2\lambda T}$ is the density for the two specified arrows and no others on $\{\circ,x\}$. The probability that this child supplies a pair is therefore
\[
p_{i,N}=\int_0^T\int_0^t \lambda^2e^{-2\lambda T}\bar q_{i+1,N}(t-s)\,ds\,dt \le \lambda^2 T^2/2 .
\]

Couple each child's contribution to a Poisson point process with intensity \eqref{eq:eligible-density}. Write $p=p_{i,N}$ and denote the Bernoulli$(p)$ count by $B$ and the Poisson$(p)$ count by $K$. Assign probability $1-p$ to $(B,K)=(0,0)$, since $\Pp(K=0)=e^{-p}\ge 1-p=\Pp(B=0)$. Assign probability $pe^{-p}$ to $(B,K)=(1,1)$, using all of $\Pp(K=1)$. On the remaining outcomes $B=1$ and $K$ is either zero or at least two. Thus
\[
\Pp(B\ne K)=1-(1-p)-pe^{-p}=p(1-e^{-p})\le p^2.
\]

When both counts are one, the conditional location law in either process is \eqref{eq:eligible-density} divided by $p$, so use the same point. Making these couplings independently for the children gives a Poisson process of intensity $n_{i,N}$ times \eqref{eq:eligible-density}, with probability at most $n_{i,N}p_{i,N}^2$ that any child's points differ.

To obtain the edges defining $\Psi_{a\sqrt{b_i}}\bar q_{i+1,N}$, the coefficient of $\bar q_{i+1,N}(t-s)\,ds\,dt$ in this intensity must change from $n_{i,N}\lambda^2e^{-2\lambda T}$ to $a^2b_i$. Use the same points at the smaller intensity and add independent Poisson points for the difference. The probability of any extra point is bounded by its expected count, at most $(T^2/2)|n_{i,N}\lambda^2e^{-2\lambda T}-a^2b_i|$.

Use the same recovery marks at $\circ$ throughout. The histories agree on $[0,T]$ unless a child edge has three or more arrows, a child's pair count differs from its Poisson count, or the intensity adjustment adds a point. Bounding these three events in that order, put
\[
\varepsilon_N(T)=\max_i\left\{ \frac{n_{i,N}(2\lambda T)^3}{6} +n_{i,N}p_{i,N}^2 +\frac{T^2}{2}\bigl|n_{i,N}\lambda^2e^{-2\lambda T}-a^2b_i\bigr| \right\}.
\]

Since $\lambda\to0$, $n_{i,N}\lambda^2=a^2b_{i,N}\to a^2b_i$ and $k$ is fixed, $\varepsilon_N(T)\to0$. The coupling gives
\begin{equation}\label{eq:approximate-fixedpoint}
\max_i\sup_{t\le T} |\bar q_{i,N}(t)-(\Psi_{a\sqrt{b_i}}\bar q_{i+1,N})(t)| \le\varepsilon_N(T).
\end{equation}

Compare \eqref{eq:approximate-fixedpoint} with \eqref{eq:cyclic-root} and put $C=a^2\max_i b_i$. For $0\le t\le T$, \eqref{eq:volterra} and its iteration give, for every $m\ge1$,
\begin{align*}
\max_i |\bar q_{i,N}(t)-q_i(t)| &\le \varepsilon_N(T)+C\int_0^t(t-u)
           \max_i |\bar q_{i,N}(u)-q_i(u)|\,du\\
&\le \varepsilon_N(T)\sum_{j=0}^{m-1}\frac{C^jt^{2j}}{(2j)!} +\frac{C^mt^{2m}}{(2m)!}.
\end{align*}

The last term uses $|\bar q_{i,N}-q_i|\le1$. Letting $m\to\infty$ gives
\[
\max_i\sup_{0\le t\le T}|\bar q_{i,N}(t)-q_i(t)| \le \varepsilon_N(T)\cosh(\sqrt C\,T)\longrightarrow0.
\]

Finally, restore the edge from $\circ$ to its parent. The processes on $S(\circ)$ and on the full tree have identical histories at $\circ$ unless an arrow occurs on that edge by $T$, an event of probability at most $2\lambda T$. This proves \eqref{eq:root-limit} for $q_{i,N}$. Likewise, deleting any $J$ edges incident to $\circ$ changes its full-tree history with probability at most $2J\lambda T\to0$.
\end{proof}

The root profiles along a fixed path determine, to leading order, the probability that infection started at one end is present at the other at time $t$.

For functions on $[0,\infty)$, write $(f*h)(t)=\int_0^t f(s)h(t-s)\,ds$.

\begin{lemma}\label{lem:geodesic}
Fix $a>0$ and an integer $m\ge0$. In each tree $\T(\bn^{(N)})$, choose a simple path $x_0,\ldots,x_m$ whose ordered vertex types $\tau_0,\ldots,\tau_m$ are the same for all $N$. Run the contact process at infection rate $\lambda_N=a\Gamma_N^{-1/(2k)}$. Under
\eqref{eq:shape-assumption}, for every finite $T>0$,
\begin{equation}\label{eq:geodesic-limit}
\lim_{N\to\infty}\sup_{0\le t\le T} \left| \lambda_N^{-m}\Pp^{\{x_0\}}(x_m\in\xi_t) -(q_{\tau_0}*\cdots*q_{\tau_m})(t) \right|=0.
\end{equation}

For $m=0$, the convolution means $q_{\tau_0}$.
\end{lemma}

\begin{proof}
For $m=0$ the claim is Lemma~\ref{lem:root-limit}. Suppose $m\ge1$. Fix $0\le t\le T$ and let $K_t$ count the arrows in both directions on the $m$ path edges up to $t$. Infection cannot reach $x_m$ from $x_0$ with fewer than $m$ arrows. To bound the contribution from more than $m$, use the fact that $K_t$ is Poisson with mean $2m\lambda_N t$. Then
\[
\lambda_N^{-m}\Pp(K_t\ge m+1) \le\lambda_N^{-m}\E\binom{K_t}{m+1} =\lambda_N\frac{(2mt)^{m+1}}{(m+1)!} \le\lambda_N\frac{(2mT)^{m+1}}{(m+1)!}.
\]

If exactly $m$ arrows occur and $x_m$ is infected at $t$, those arrows must be $x_{j-1}\to x_j$ at times $0<s_1<\cdots<s_m<t$. The density of these arrows with no other arrows on the path is $\lambda_N^m e^{-2m\lambda_N t}\,ds_1\cdots ds_m$. Cut the path edges and write $C_j$ for the component containing $x_j$. Given the arrows, infection reaches $(x_m,t)$ exactly when, in each $C_j$, there is an infection path from $(x_j,s_j)$ to $(x_j,s_{j+1})$, where $s_0=0$ and $s_{m+1}=t$. The $m+1$ events are independent. The probability of the $j$th event is $f_{j,N}(s_{j+1}-s_j)=\Pp_{C_j}^{\{x_j\}}(x_j\in\xi_{s_{j+1}-s_j})$.

Integrating over the arrow times gives the first equality below. The second follows by expanding the definition of convolution, with $s_m$ ranging from $0$ to $t$, $s_{m-1}$ from $0$ to $s_m$, and so on.
\begin{align*}
\lambda_N^{-m}\Pp^{\{x_0\}}(x_m\in\xi_t,\,K_t=m) &=\int_{0<s_1<\cdots<s_m<t} e^{-2m\lambda_N t}
          \prod_{j=0}^m f_{j,N}(s_{j+1}-s_j)\,ds_1\cdots ds_m,\\
(q_{\tau_0}*\cdots*q_{\tau_m})(t) &=\int_{0<s_1<\cdots<s_m<t} \prod_{j=0}^m q_{\tau_j}(s_{j+1}-s_j)\,ds_1\cdots ds_m .
\end{align*}

At most two edges incident to each $x_j$ were cut. By Lemma~\ref{lem:root-limit}, there are numbers $\delta_N\to0$ such that $|f_{j,N}(u)-q_{\tau_j}(u)|\le\delta_N$ for every $j$ and $0\le u\le T$. Fix the crossing times and write $u_j=s_{j+1}-s_j$. Since all factors lie in $[0,1]$, replacing them one at a time gives
\[
\left|\prod_{j=0}^m f_{j,N}(u_j)-\prod_{j=0}^m q_{\tau_j}(u_j)\right|
\le \sum_{j=0}^m |f_{j,N}(u_j)-q_{\tau_j}(u_j)|
\le (m+1)\delta_N .
\]

The factor $e^{-2m\lambda_N t}$ changes the first product by at most $1-e^{-2m\lambda_N t}\le 2m\lambda_N T$. The ordered crossing times lie in $[0,T]^m$, so for every $0\le t\le T$,
\begin{align*}
&\left|\lambda_N^{-m}\Pp^{\{x_0\}}(x_m\in\xi_t)
       -(q_{\tau_0}*\cdots*q_{\tau_m})(t)\right|\\
&\quad\le \lambda_N^{-m}\Pp(K_t\ge m+1)
       +\int_{[0,T]^m}\bigl(2m\lambda_N T+(m+1)\delta_N\bigr)\,ds_1\cdots ds_m\\
&\quad\le \lambda_N\frac{(2mT)^{m+1}}{(m+1)!}
       +T^m\bigl(2m\lambda_N T+(m+1)\delta_N\bigr)
       \xrightarrow[N\to\infty]{}0. \qedhere
\end{align*}
\end{proof}

Now we set up the supermartingale weights, as we would for the branching random walk upper bound. Fix $\lambda=a\Gamma_N^{-1/(2k)}$. For each type $i$, choose a starting vertex $\circ$ of that type in $\T(\bn^{(N)})$. Give $\circ$ weight $w_{\circ}=1$ and set
\begin{equation}\label{eq:weights}
\frac{w_y}{w_x}=n_{j,N}^{-1/2} \quad\text{if $x$ has type $j$ and $y$ is a child of $x$.}
\end{equation}

These ratios specify weights uniquely on the tree. At a type-$j$ vertex $x$,
\[
\sum_{y:y\sim x}w_y =(\sqrt{n_{j-1,N}}+\sqrt{n_{j,N}})w_x.
\]

Put
\begin{equation}\label{eq:V}
V_N={\lambda}\max_i (\sqrt{n_{i-1,N}}+\sqrt{n_{i,N}}) \longrightarrow aB,\qquad B=\max_i(\sqrt{b_{i-1}}+\sqrt{b_i}).
\end{equation}

For the contact process started from $\{\circ\}$, define its expected weighted infected mass at time $t\ge0$ by
\[
Z_{i,N}(t):=\E^{\{\circ\}}\!\left[\sum_{x\in\xi_t} w_x\right] =\sum_{x\in V(\T(\bn^{(N)}))} w_x\,\Pp^{\{\circ\}}(x\in\xi_t).
\]

An automorphism preserving types and the distinguished end carries $\circ$ to any other type-$i$ vertex and transports the normalised weights to those for that vertex. Thus the definition is independent of the choice of $\circ$.

\begin{lemma}\label{lem:mass-limit}
For each starting type $i$ there is a limiting expected weighted infected mass $Z_i(t)$ such that, for every finite $T>0$,
\begin{equation}\label{eq:mass-limit}
\lim_{N\to\infty}\sup_{0\le t\le T}|Z_{i,N}(t)-Z_i(t)|=0.
\end{equation}

The limit satisfies $q_i(t)\le Z_i(t)\le e^{(aB-1)t}$. For fixed $\bb$ and $t$, $Z_i(t)$ is continuous in $a$.
\end{lemma}

\begin{proof}
Group the vertices in the sum defining $Z_{i,N}(t)$ by the numbers of parent and child steps on the geodesic from $\circ$. Once the path steps to a child, a step to its parent would retrace the last edge. Thus the $\ell$ parent steps come first, followed by the $r$ child steps. Let $V_{\ell,r}$ be the set of endpoints of these paths and put $m=\ell+r$. All paths in the class have the same vertex types $\tau_0,\ldots,\tau_m$, starting at $i$ and decreasing by one on parent steps and increasing by one on child steps, mod $k$. Fix $\ell,r$. For all sufficiently large $N$, choose $y_N\in V_{\ell,r}$. Lemma~\ref{lem:geodesic} states that
\[
\lim_{N\to\infty}\sup_{0\le t\le T} \left|\frac{\Pp^{\{\circ\}}(y_N\in\xi_t)}{\lambda^m} -(q_{\tau_0}*\cdots*q_{\tau_m})(t)\right|=0.
\]

There is one profile for each vertex, from the time of arrival to the next path arrow, or to time $t$ at the last vertex.

For the $j$th edge let $p_j$ be the type of its parent vertex. A step to that parent is forced and multiplies the weight by $\sqrt{n_{p_j,N}}$. In the other direction there are $n_{p_j,N}$ children of relative weight $1/\sqrt{n_{p_j,N}}$. If $\ell,r>0$, the first child step has just $n_{i-\ell,N}-1$ choices, since it must avoid the edge just traversed. Using $\lambda\sqrt{n_{p_j,N}}=a\sqrt{b_{p_j,N}}$, we obtain
\[
\lambda^m\sum_{x\in V_{\ell,r}}w_x =\left(1-\frac{\1_{\{\ell>0,\,r>0\}}}{n_{i-\ell,N}}\right) \prod_{j=1}^{m}a\sqrt{b_{p_j,N}} \xrightarrow[N\to\infty]{}\prod_{j=1}^{m}a\sqrt{b_{p_j}}.
\]

By symmetry every $x\in V_{\ell,r}$ has the same infection probability as $y_N$. Combining the two limits gives, uniformly for $0\le t\le T$,
\begin{align*}
\sum_{x\in V_{\ell,r}}w_x\,\Pp^{\{\circ\}}(x\in\xi_t) &=\left(\lambda^m\sum_{x\in V_{\ell,r}}w_x\right)
   \frac{\Pp^{\{\circ\}}(y_N\in\xi_t)}{\lambda^m}\\
&\xrightarrow[N\to\infty]{} \left(\prod_{j=1}^{m}a\sqrt{b_{p_j}}\right) (q_{\tau_0}*\cdots*q_{\tau_m})(t).
\end{align*}

For $\ell$ parent steps followed by $r$ child steps, let $m=\ell+r$ and write $s_j$ for the time the path crosses its $j$th edge, with $s_0=0$ and $s_{m+1}=t$. Writing $q_j(u;a)$ to show the dependence on the rate, we define the limiting mass by
\begin{equation}\label{eq:weighted-path-sum}
Z_i(t)=\sum_{\ell,r\ge0} a^m\left(\prod_{j=1}^m\sqrt{b_{p_j}}\right) \int_{0<s_1<\cdots<s_m<t} \prod_{j=0}^m q_{\tau_j}(s_{j+1}-s_j;a)\,ds_1\cdots ds_m .
\end{equation}

For $\ell=r=0$ there are no crossings and the summand is $q_i(t;a)$.

For $0\le t\le T$, Lemma~\ref{lem:weighted-brw} bounds the tail by
\begin{equation}\label{eq:uniform-tail}
\sum_{x:\,d(\circ,x)\ge K}w_x\,\Pp^{\{\circ\}}(x\in\xi_t) \le e^{-t}\sum_{j\ge K}\frac{(V_Nt)^j}{j!} \le\sum_{j\ge K}\frac{(T\sup_N V_N)^j}{j!}.
\end{equation}

The same bound holds for the sum in \eqref{eq:weighted-path-sum} over $K\le\ell+r\le K'$ by convergence of each of these finitely many terms. These partial sums increase to the full tail $\ell+r\ge K$ as $K'\to\infty$. The terms with $\ell+r<K$ converge uniformly, and $\sup_N V_N<\infty$ by \eqref{eq:V}. Thus
\[
\limsup_{N\to\infty}\sup_{0\le t\le T}|Z_{i,N}(t)-Z_i(t)| \le 2\sum_{j\ge K}\frac{(T\sup_N V_N)^j}{j!} \xrightarrow[K\to\infty]{}0.
\]

The root term and the total-mass bound give
\[
q_i(t)\le Z_i(t)=\lim_{N\to\infty}Z_{i,N}(t) \le\lim_{N\to\infty}e^{(V_N-1)t}=e^{(aB-1)t}.
\]

Fix $t$ and $\bb$. For fixed crossing times in \eqref{eq:weighted-path-sum}, each root probability is continuous in $a$ by Lemma~\ref{lem:fixed-point}. Their product is bounded by $1$ on a subset of $[0,t]^m$, so each integral is continuous in $a$ by bounded convergence. Its prefactor is a constant times $a^m$, so each summand is continuous as well.

By \eqref{eq:V}, $V_N$ is proportional to $a$, so the right-hand side of \eqref{eq:uniform-tail} at $a=A$ bounds the tail of \eqref{eq:weighted-path-sum} for every $0<a\le A$. Thus the series converges uniformly on $(0,A]$, and its sum $Z_i(t)$ is continuous in $a$.
\end{proof}

Finally we use these branching random walk calculations for a contraction test to give local extinction.

\begin{proposition}\label{prop:extinction-test}
If $H_{\bb}(a)<1$, then the contact process at ${\lambda}=a\Gamma_N^{-1/(2k)}$ dies out locally for all sufficiently large $N$.
\end{proposition}

\begin{proof}
Write $H=H_{\bb}(a)$ and $\chi_i=\chi_i(a,\bb)$. Since $H<1$ and every $\chi_i\ge1$, all susceptibilities are finite. For a path of length $m$, let $u_0,\ldots,u_m$ be the lengths of the intervals into which the $m$ crossing times divide $[0,t]$. The change from the crossing times and $t$ to these lengths has Jacobian $1$, so Tonelli's theorem gives
\begin{align*}
\int_0^\infty(q_{\tau_0}*\cdots*q_{\tau_m})(t)\,dt
 &=\int_{[0,\infty)^{m+1}}\prod_{j=0}^m q_{\tau_j}(u_j)\,du_0\cdots du_m\\
&=\prod_{j=0}^m\chi_{\tau_j}.
\end{align*}

In \eqref{eq:weighted-path-sum} the starting vertex contributes $\chi_i$. A step to a parent of type $j$ contributes $a\sqrt{b_j}\chi_j$, and a step from that parent to a child contributes $a\sqrt{b_j}\chi_{j+1}$. Consequently the integral of the nonnegative series is
\[
\int_0^\infty Z_i(t)\,dt =\chi_i\sum_{\ell,r\ge0} \left(\prod_{j=1}^{\ell}a\sqrt{b_{i-j}}\,\chi_{i-j}\right) \left(\prod_{j=0}^{r-1}a\sqrt{b_{i-\ell+j}}\,\chi_{i-\ell+j+1}\right).
\]

For a full period of steps in either direction the product is
\[
\prod_{j=1}^k\bigl(a\sqrt{b_j}\chi_j\bigr) =\prod_{j=1}^k\bigl(a\sqrt{b_j}\chi_{j+1}\bigr) =a^k\prod_{j=1}^k\chi_j=H,
\]
since the types are cyclic and $\prod_{j=1}^k b_j=1$. Write $\ell=uk+s$ and $r=vk+h$ with $0\le s,h<k$. The child steps start at type $i-\ell\equiv i-s\pmod{k}$. Grouping the full periods and then summing the two geometric series gives
\begin{align}
\int_0^\infty Z_i(t)\,dt &=\chi_i\sum_{u,v\ge0}H^{u+v} \sum_{s,h=0}^{k-1} \left(\prod_{j=1}^s a\sqrt{b_{i-j}}\,\chi_{i-j}\right)
   \left(\prod_{j=0}^{h-1}a\sqrt{b_{i-s+j}}\,\chi_{i-s+j+1}\right)\notag\\
&=\frac{\chi_i}{(1-H)^2} \sum_{s,h=0}^{k-1} \left(\prod_{j=1}^s a\sqrt{b_{i-j}}\,\chi_{i-j}\right) \left(\prod_{j=0}^{h-1}a\sqrt{b_{i-s+j}}\,\chi_{i-s+j+1}\right) <\infty .
 \label{eq:integrated-mass}
\end{align}

Fix $N$ and $T>0$, put $\delta_N=\max_i Z_{i,N}(T)$, and use the weights \eqref{eq:weights} throughout the tree. An automorphism preserving types and the distinguished end carries a process from a type-$i$ representative to one from any type-$i$ vertex $x$, multiplying every weight by $w_x$. Additivity and a union bound therefore give, for every finite $A$,
\begin{equation}\label{eq:contraction}
\E^A\sum_v w_v\1\{v\in\xi_T\} \le\sum_{x\in A}w_x Z_{\operatorname{type}(x),N}(T) \le\delta_N\sum_{x\in A}w_x.
\end{equation}

As in Huang's proof of Theorem~1.7 \cite[p.~19]{MR4125782}, condition on $\xi_s$ and apply \eqref{eq:contraction} for the next $T$ units of time. By the Markov property,
\[
Z_{i,N}(s+T)
\le \E^{\{\circ\}}\!\left[\sum_{x\in\xi_s}w_x Z_{\operatorname{type}(x),N}(T)\right]
\le \delta_N Z_{i,N}(s).
\]

Starting at $s=u$ and using Lemma~\ref{lem:weighted-brw} to bound $Z_{i,N}(u)$ gives, for $j\ge0$ and $0\le u<T$,
\begin{equation}\label{eq:iterated-contraction}
Z_{i,N}(jT+u)\le\delta_N^j Z_{i,N}(u)\le\delta_N^j e^{(V_N-1)u}.
\end{equation}

By \eqref{eq:integrated-mass}, the finite sum $\sum_i Z_i$ is integrable, so there is a $T>0$ with $\sum_i Z_i(T)<1$. At this $T$, Lemma~\ref{lem:mass-limit} gives
\[
\lim_{N\to\infty}\delta_N=\max_i Z_i(T)<1.
\]

For all sufficiently large $N$, $\delta_N<1$. Since the root has weight $1$, \eqref{eq:iterated-contraction} gives
\[
\int_0^\infty q_{i,N}(t)\,dt \le \int_0^\infty Z_{i,N}(t)\,dt
\le \sum_{j=0}^\infty\delta_N^j\int_0^T e^{(V_N-1)u}\,du<\infty .
\]

Equation~\eqref{eq:occupation-extinction} proves local extinction.
\end{proof}


\section{Characterisation of the limiting constant}
\label{sec:threshold}

Choose a root $e_0=\circ$ of type $i$ and a forward ray $e_0,e_1,e_2,\ldots$, and let $S(x)$ be the subtree consisting of $x$ and all its descendants. We briefly give an argument from \cite{MR1717346}. Define
\[
R_m=\Pp^{\{\circ\}}(e_{km}\text{ is ever infected}), \qquad V_m=\sup_{t\ge0} \Pp_{S(\circ)}^{\{\circ\}}(e_{km}\in\xi_t).
\]

By the strong Markov property at the first hit of $e_{km}$, deleting all other infections at that time,
\[
R_{m+j}\ge R_mR_j.
\]

Vertices separated by a multiple of $k$ on the ray have the same rooted, end-oriented full-tree structure. The supermultiplicative limit thus exists and satisfies
\begin{equation}\label{eq:beta}
\beta=\lim_{m\to\infty}R_m^{1/m} =\sup_{m\ge1}R_m^{1/m}.
\end{equation}

Let $h_i(\lambda)$ be the full-tree probability that a process started at one type-$i$ vertex ever infects a prescribed forward child. Child symmetry makes the choice of that child immaterial.

\begin{lemma}\label{lem:hitting-criterion}
On $\T(\bn)$ with arbitrary $\lambda>0$ and $\Gamma=\prod_i n_i$,
\begin{equation}\label{eq:hitting-criterion}
\prod_{i=1}^k h_i(\lambda)>\Gamma^{-1/2} \quad\Longrightarrow\quad \text{local survival}.
\end{equation}
\end{lemma}

\begin{proof}
Successively restart the process with only the next vertex on the forward ray infected. The strong Markov property gives $R_1\ge\prod_i h_i(\lambda)$, so
\eqref{eq:beta} implies $\beta>\Gamma^{-1/2}$.

In the notation immediately before this lemma, \cite[Lemma~6.1]{MR4125782} states that
\[
\lim_{m\to\infty}V_m^{1/m}=\beta.
\]

Since $\beta>\Gamma^{-1/2}$, we can choose $L\ge1$ with $V_L>\Gamma^{-L/2}$. By the definition of $V_L$ as a supremum over times, some finite $t>0$ then satisfies
\begin{equation}\label{eq:block-choice}
p:=\Pp^{\{\circ\}}_{S(\circ)}(e_{kL}\in\xi_t)>\Gamma^{-L/2}, \qquad\text{so}\qquad \Gamma^L p^2>1.
\end{equation}

At time $t$ each of the $\Gamma^L$ descendants at depth $kL$ is infected with probability $p$. Retaining only these infections and restarting in their disjoint forward subtrees gives a branching process with mean offspring number $\Gamma^L p$. After $j$ blocks its mean size is $(\Gamma^L p)^j$. For a specified descendant at depth $jkL$, concatenating $j$ blocks along its ancestral line gives an infection probability of at least $p^j$ at time $jt$. By self-duality on $S(\circ)$, a process started at that descendant also infects $\circ$ at time $jt$ with probability at least $p^j$. Thus the product of the mean size and the return probability is at least
\[
(\Gamma^L p)^j p^j=(\Gamma^L p^2)^j\longrightarrow\infty.
\]

This is the condition used in Huang's branching-and-return argument \cite[Section~6.3, (6.10)--(6.12)]{MR4125782}, which gives local survival in $S(\circ)$ and hence on the full tree.
\end{proof}

This result provides the following result, the supercritical susceptibility test.

\begin{proposition}\label{prop:survival-test}
If $H_{\bb}(a)>1$ then the contact process at ${\lambda}=a\Gamma_N^{-1/(2k)}$ survives locally for all sufficiently large $N$.
\end{proposition}

\begin{proof}
Choose finite $T$ such that
\begin{equation}\label{eq:finite-occupation-test}
a^k\prod_{i=1}^k\int_0^Tq_i(t)\,dt>1.
\end{equation}

For each type-$i$ root, reserve one forward child as a target. Let $C$ be the root component after removing that edge, and run the process in $C$ from the root alone. Let $A_{i,N}(T)$ be the time the root spends infected in this component up to time $T$. Then $0\le A_{i,N}(T)\le T$ and, by Lemma~\ref{lem:root-limit},
\[
\E A_{i,N}(T)=\int_0^T \Pp_C^{\{\circ\}}(\circ\in\xi_t)\,dt
\longrightarrow\int_0^Tq_i(t)\,dt.
\]

Restore only the independent arrows from the root to the target, suppressing all feedback. Conditional on the root path, the probability of at least one such transmission by time $T$ is $1-e^{-{\lambda} A_{i,N}(T)}$. The full-tree hitting probability is therefore at least
\[
h_i({\lambda}) \ge\E[1-e^{-{\lambda} A_{i,N}(T)}] \ge{\lambda}\E A_{i,N}(T)-{\lambda}^2T^2/2.
\]

Since $\lambda=a\Gamma_N^{-1/(2k)}\to0$, for each type $i$ we have
\[
\liminf_{N\to\infty}\frac{h_i(\lambda)}{\lambda} \ge \lim_{N\to\infty} \left(\E A_{i,N}(T)-\frac{\lambda T^2}{2}\right) =\int_0^T q_i(t)\,dt>0.
\]

Using $\sqrt{\Gamma_N}\lambda^k=a^k$ and multiplying over the $k$ types gives
\[
\liminf_{N\to\infty} \sqrt{\Gamma_N}\prod_{i=1}^k h_i({\lambda}) \ge a^k\prod_{i=1}^k\int_0^Tq_i(t)\,dt>1.
\]

Lemma~\ref{lem:hitting-criterion} proves the claim.
\end{proof}

\begin{proof}[Proof of the characterisation in Theorem~\ref{thm:periodic}]
Fix $\bb$. For $0<a_1<a_2$, monotonicity of the susceptibilities gives
\begin{equation}\label{eq:strict-crossing}
H_{\bb}(a_2)\ge(a_2/a_1)^kH_{\bb}(a_1).
\end{equation}

If $H_{\bb}(\kappa)=1$, then
\[
H_{\bb}(a)\le(a/\kappa)^k<1 \quad(0<a<\kappa),\qquad H_{\bb}(a)\ge(a/\kappa)^k>1 \quad(a>\kappa).
\]

Thus \eqref{eq:kappa} has at most one solution.

For existence, we show that $\{a>0:H_{\bb}(a)>1\}$ and $\{a>0:H_{\bb}(a)<1\}$ are open and nonempty. Write $q_i(t;a)$ and $Z_i(t;a)$ for the root profile and the limiting weighted mass at rate $a$. If $H_{\bb}(a)>1$, choose a finite $T$ for which \eqref{eq:finite-occupation-test} holds. The profiles lie in $[0,1]$ and are continuous in $a$ at each time, by Lemma~\ref{lem:fixed-point}. By bounded convergence on $[0,T]$ the product is continuous in $a$, so $\{a:H_{\bb}(a)>1\}$ is open.

For the other set we claim that
\[
H_{\bb}(a)<1 \quad\Longleftrightarrow\quad \max_i Z_i(T;a)<1 \text{ for some } T>0.
\]

If $H_{\bb}(a)<1$, then $\int_0^\infty\sum_i Z_i(t;a)\,dt<\infty$ by \eqref{eq:integrated-mass}, so such a $T$ exists.

Now suppose that $\delta:=\max_i Z_i(T;a)<1$ for some $T>0$, to show that $H_{\bb}(a)<1$. For $j\ge0$ and $0\le u<T$, \eqref{eq:iterated-contraction},
\eqref{eq:mass-limit} and \eqref{eq:V} give
\[
Z_i(jT+u;a)=\lim_{N\to\infty}Z_{i,N}(jT+u;a) \le \delta^j e^{(aB-1)u}.
\]

In particular,
\[
\int_0^\infty Z_i(t;a)\,dt \le \sum_{j=0}^\infty\delta^j\int_0^T e^{(aB-1)u}\,du<\infty .
\]

For a lower bound, keep only $(\ell,r)=(0,mk)$ in \eqref{eq:weighted-path-sum}. With $u_j=s_{j+1}-s_j$, we have
\begin{align*}
\int_0^\infty Z_i(t;a)\,dt &\ge \sum_{m=0}^\infty \left(\prod_{j=0}^{mk-1}a\sqrt{b_{i+j}}\right)
   \prod_{j=0}^{mk}\left(\int_0^\infty q_{i+j}(u_j;a)\,du_j\right)\\
&=\chi_i(a,\bb)\sum_{m=0}^\infty \left(a^k\prod_{j=1}^k\sqrt{b_j}\,\chi_j(a,\bb)\right)^m =\chi_i(a,\bb)\sum_{m=0}^\infty H_{\bb}(a)^m ,
\end{align*}
since the types repeat every $k$ steps and $\prod_{j=1}^k b_j=1$. Since $\chi_i(a,\bb)\ge1$, this forces $H_{\bb}(a)<1$. For each fixed $T$ the maximum $\max_i Z_i(T;a)$ is continuous in $a$ by Lemma~\ref{lem:mass-limit}. Hence $\{a:H_{\bb}(a)<1\}$ is open as well.

Both strict inequalities occur. For $0<a<1/(1+B)$, Lemma~\ref{lem:mass-limit} gives
\[
\chi_i(a,\bb)\le\int_0^\infty e^{-(1-aB)t}\,dt=\frac{1}{1-aB}, \qquad H_{\bb}(a)\le\left(\frac{a}{1-aB}\right)^k<1.
\]

For $a>1$ we have $H_{\bb}(a)\ge a^k>1$ since every $\chi_i\ge1$. The interval $(0,\infty)$ cannot be the union of these two disjoint nonempty open sets. Hence $H_{\bb}(\kappa)=1$ for some $\kappa>0$.
\end{proof}

Applying the extinction and survival tests on either side of this solution gives the convergence in Theorem~\ref{thm:periodic}.

\begin{proof}[Proof of convergence in Theorem~\ref{thm:periodic}]
Let $\kappa=\kappa(\bb)$ be the unique solution of \eqref{eq:kappa} and fix $0<\varepsilon<\kappa$. Since
\[
H_{\bb}(\kappa-\varepsilon)<1<H_{\bb}(\kappa+\varepsilon),
\]
Proposition~\ref{prop:extinction-test} at $\kappa-\varepsilon$ and Proposition~\ref{prop:survival-test} at $\kappa+\varepsilon$ give, for all sufficiently large $N$,
\[
\kappa-\varepsilon \le \Gamma_N^{1/(2k)}\lambda_2(\T(\bn^{(N)})) \le \kappa+\varepsilon .
\]

This proves \eqref{eq:periodic-limit}. The same limit holds for the normalisation by degrees, since
\[
\frac{(d_{1,N}\cdots d_{k,N})^{1/(2k)}}{\Gamma_N^{1/(2k)}} =\prod_{i=1}^k\left(1+\frac{1}{n_{i,N}}\right)^{1/(2k)} \longrightarrow 1. \qedhere
\]
\end{proof}

We now use the susceptibility characterisation to prove the two bounds in Proposition~\ref{prop:upper-bounds}. For the lower bound, an edge can start a new infected period only if a recovery occurs between its two endpoints.

\begin{lemma}\label{lem:effective-return}
If $c^2\int_0^\infty f(u)(1-e^{-u})\,du<1$, then
\[
\int_0^\infty \Psi_c f(t)\,dt \le \frac{1}{1-c^2\int_0^\infty f(u)(1-e^{-u})\,du}.
\]
\end{lemma}

\begin{proof}
Reveal each recovery as it occurs and each edge $(v,w)$ at its later time $w$. At time $w$, edges with starting times in $dv$ arrive at rate $c^2f(w-v)\,dv$. Such an edge starts a new infected period exactly when the root was infected at $v$ and is healthy just before $w$. Counting every edge from an infected root that spans a recovery also counts returns for which the root has already been reinfected. Since recovery marks after $v$ are independent of the state at $v$, the mean rate $r(w)$ at which new infected periods start at $w>0$ satisfies
\begin{align*}
 r(w)&=c^2\int_0^w f(w-v)\,\Pp\bigl(X^{c,f}(v)=1,\,X^{c,f}(w-)=0\bigr)\,dv\\
&\le c^2\int_0^w f(w-v)\Psi_c f(v)(1-e^{-(w-v)})\,dv .
\end{align*}

An infected period starting at $w$ lasts through $t$ if there is no recovery in $(w,t]$, an independent event of probability $e^{-(t-w)}$. Including the initial infected period gives
\[
\Psi_c f(t)=e^{-t}+\int_0^t r(w)e^{-(t-w)}\,dw .
\]

For $T>0$, the bound on $r(w)$ and the substitution $u=w-v$ give
\begin{align*}
\int_0^T \Psi_c f(t)\,dt
 &=1-e^{-T}+\int_0^T r(w)\bigl(1-e^{-(T-w)}\bigr)\,dw\\
 &\le 1+\int_0^T r(w)\,dw\\
&\le 1+c^2\int_0^T \Psi_c f(v)
       \int_0^{T-v} f(u)(1-e^{-u})\,du\,dv\\
&\le 1+c^2\left(\int_0^T \Psi_c f(v)\,dv\right) \left(\int_0^\infty f(u)(1-e^{-u})\,du\right).
\end{align*}

Since $c^2\int_0^\infty f(u)(1-e^{-u})\,du<1$,
\[
\int_0^\infty \Psi_c f(t)\,dt \le \frac{1}{1-c^2\int_0^\infty f(u)(1-e^{-u})\,du}. \qedhere
\]
\end{proof}

Iterating this estimate around the period gives upper bounds on the root susceptibilities and hence a lower bound on $\kappa(\bb)$.

\begin{proof}[Proof of the lower bound in Proposition~\ref{prop:upper-bounds}]
Take $a=(4-2\sqrt2)/\rho(M_{\bb})$. We want to show that $H_{\bb}(a)=a^k\prod_i\chi_i<1$. The lower iterates in Lemma~\ref{lem:fixed-point} are given by
\[
f_i^{(0)}(t)=e^{-t},\qquad f_{i-1}^{(j+1)}=\Psi_{a\sqrt{b_{i-1}}} f_i^{(j)} \quad(j\ge0),
\]
with type indices read mod $k$. We shall find numbers $x_i\ge1$ with $\prod_i x_i=a^{-k}$ such that every iterate $f_i^{(j)}$ has integral at most $x_i$.

It suffices to choose the $x_i$ so that, for each $i$,
\begin{equation}\label{eq:lower-bound-iteration}
\frac1{x_{i-1}}+a^2b_{i-1}(x_i-\tfrac12)\le1.
\end{equation}

For the induction step, suppose $\int_0^\infty f_i^{(j)}(t)\,dt\le x_i$. Since $f_i^{(j)}(t)\ge e^{-t}$,
\[
\int_0^\infty f_i^{(j)}(t)(1-e^{-t})\,dt =\int_0^\infty f_i^{(j)}(t)\,dt-\int_0^\infty f_i^{(j)}(t)e^{-t}\,dt \le x_i-\int_0^\infty e^{-2t}\,dt=x_i-\tfrac12.
\]

By \eqref{eq:lower-bound-iteration}, $1-a^2b_{i-1}(x_i-\tfrac12)\ge 1/x_{i-1}>0$. Lemma~\ref{lem:effective-return} therefore gives
\begin{align*}
\int_0^\infty f_{i-1}^{(j+1)}(t)\,dt
 &=\int_0^\infty (\Psi_{a\sqrt{b_{i-1}}}f_i^{(j)})(t)\,dt\\
&\le \frac{1}{1-a^2b_{i-1}\int_0^\infty f_i^{(j)}(u)(1-e^{-u})\,du} \le \frac{1}{1-a^2b_{i-1}(x_i-\tfrac12)} \le x_{i-1}.
\end{align*}

To construct the $x_i$, divide row $i$ of $aM_{\bb}$ by $1+a^2b_i/2$ for each $i$. The resulting matrix is nonnegative and irreducible, since every type is connected to the next. Take its positive Perron eigenvector $v$, normalised by $\sum_i v_i^2=1$, and let $\mu$ be its eigenvalue. Multiplying $a(M_{\bb}v)_i=\mu(1+a^2b_i/2)v_i$ by $v_i$ and summing gives
\[
u:=a v^{\mathsf T}M_{\bb}v =\mu\left(1+\tfrac{a^2}{2}\sum_i b_i v_i^2\right).
\]

Since $M_{\bb}$ is symmetric and $v$ is a unit vector,
\[
0\le u=2a\sum_i\sqrt{b_i}\,v_i v_{i+1} \le a\rho(M_{\bb})=4-2\sqrt2.
\]

By Cauchy--Schwarz and $\sum_i v_{i+1}^2=1$,
\begin{align*}
1+\frac{a^2}{2}\sum_i b_i v_i^2 &\ge 1+\frac{a^2}{2}\left(\sum_i\sqrt{b_i}\,v_i v_{i+1}\right)^2
 =1+\frac{u^2}{8}\\
&=u+\frac{(4-2\sqrt2-u)(4+2\sqrt2-u)}{8}\ge u.
\end{align*}

Thus
\[
\mu=\frac{u}{1+\frac{a^2}{2}\sum_i b_i v_i^2}\le 1.
\]

Set $x_i=v_i/(v_{i-1}a\sqrt{b_{i-1}})$, and note $\prod_i x_i=a^{-k}$ since $\prod_i b_i=1$. The eigenvector equation in row $i-1$ now gives
\begin{align*}
\frac1{x_{i-1}}+a^2b_{i-1}x_i &=\frac{a\sqrt{b_{i-2}}v_{i-2}+a\sqrt{b_{i-1}}v_i}{v_{i-1}}
 =\frac{a(M_{\bb}v)_{i-1}}{v_{i-1}}\\
&=\mu(1+a^2b_{i-1}/2)\le 1+a^2b_{i-1}/2,
\end{align*}
which proves \eqref{eq:lower-bound-iteration}.

The induction step carries the bounds $\int_0^\infty f_i^{(j)}(t)\,dt\le x_i$ from $j$ to $j+1$. For the $x_i$ just constructed, we still need the base case
\[
\int_0^\infty f_i^{(0)}(t)\,dt\le x_i \qquad (1\le i\le k).
\]

The left-hand side is $\int_0^\infty e^{-t}\,dt=1$ for every type, so it remains to prove $1\le x_i$ for each $i$. Since $M_{\bb}$ is symmetric, its largest eigenvalue satisfies
\begin{align*}
\rho(M_{\bb}) &=\max_{z\ne0}\frac{z^{\rm T}M_{\bb}z}{z^{\rm T}z}
 \ge\frac{\1^{\rm T}M_{\bb}\1}{\1^{\rm T}\1}\\
&=\frac1k\sum_{i=1}^k(\sqrt{b_{i-1}}+\sqrt{b_i}) =\frac2k\sum_{i=1}^k\sqrt{b_i} \ge 2\left(\prod_{i=1}^k b_i\right)^{1/(2k)}=2,
\end{align*}
where the last inequality is the arithmetic--geometric mean inequality. Thus $a<1$ and $\prod_i x_i=a^{-k}>1$. If $x_{i-1}\le1$ for some $i$, then $1/x_{i-1}\ge1$. Since $a^2b_{i-1}>0$, \eqref{eq:lower-bound-iteration} gives
\[
x_i\le \frac12+\frac{1-1/x_{i-1}}{a^2b_{i-1}}\le\frac12.
\]

Repeating this around the period would give $x_j\le1$ for every $j$, contradicting the product. Thus every $x_i>1$ and the induction applies. For each $t\ge0$, the nonnegative iterates $f_i^{(j)}(t)$ increase to $q_i(t)$ as $j\to\infty$. By the monotone convergence theorem,
\[
\chi_i=\int_0^\infty q_i(t)\,dt =\lim_{j\to\infty}\int_0^\infty f_i^{(j)}(t)\,dt\le x_i.
\]

Fix $t>0$. The root is infected at $t$ if there is no recovery in $(0,t)$. If there is exactly one recovery and it lies in the middle third of $(0,t)$, any edge from the first third to the last carries infection across it. This recovery pattern has probability $(t/3)e^{-t}$, independently of the edges. Since $q_{i+1}(w-v)\ge e^{-(w-v)}\ge e^{-t}$ for $0<v<w<t$,
\begin{align*}
q_i(t) &\ge e^{-t}+\frac{t}{3}e^{-t} \left[1-\exp\left(-a^2b_i\int_0^{t/3}\int_{2t/3}^t
                  q_{i+1}(w-v)\,dw\,dv\right)\right]\\
&\ge e^{-t}+\frac{t}{3}e^{-t} \left[1-\exp\left(-a^2b_i\,\frac{t^2}{9}e^{-t}\right)\right]>e^{-t}.
\end{align*}

Thus $\int_0^\infty e^{-t}q_i(t)\,dt>\int_0^\infty e^{-2t}\,dt=\tfrac12$. Applying Lemma~\ref{lem:effective-return} to $q_{i-1}=\Psi_{a\sqrt{b_{i-1}}}q_i$ and using \eqref{eq:lower-bound-iteration} gives
\begin{align*}
\frac1{\chi_{i-1}}
 &\ge 1-a^2b_{i-1}\left(\chi_i-\int_0^\infty e^{-t}q_i(t)\,dt\right)\\
&>1-a^2b_{i-1}(x_i-\tfrac12)\ge\frac1{x_{i-1}}.
\end{align*}

Thus \[H_{\bb}(a)=a^k\prod_i\chi_i<a^k\prod_i x_i=1=H_{\bb}(\kappa(\bb)).\]

Since $H_{\bb}$ is nondecreasing, $a<\kappa(\bb)$.
\end{proof}

For the upper bound we count edges that can still return to the root. Fix $a>0$ and write $q_i$ and $\chi_i$ for the root profiles and susceptibilities in \eqref{eq:cyclic-root} and \eqref{eq:H}. If any $\chi_i$ is infinite, then $H_{\bb}(a)>1$ already. We therefore consider the case where all susceptibilities are finite.

\begin{lemma}\label{lem:pending-returns}
Let $\alpha_i=a^2b_i\chi_{i+1}$. For each type $i$,
\begin{equation}\label{eq:pending-levels}
\chi_i\ge 1+ \sum_{j=1}^\infty \frac{\alpha_i^j}{\prod_{\ell=1}^j(1+\ell/\chi_{i+1})}.
\end{equation}
\end{lemma}

\begin{proof}
For every type $i$, $q_i(u+t)\ge q_i(u)q_i(t)$. Indeed, if the root is infected at $u$, we can restart there using only edges starting after $u$ and recovery marks after $u$. These are independent of the root state at $u$. Integration in $t$ gives
\[
\int_u^\infty q_i(v)\,dv\ge\chi_i q_i(u).
\]

With $s$ the start time and $v=w-s$ the delay, the edge intensity in \eqref{eq:edge-intensity} for a type-$i$ root factors as
\[
a^2 b_i q_{i+1}(v)\,ds\,dv
=\alpha_i\,ds\,\frac{q_{i+1}(v)}{\chi_{i+1}}\,dv.
\]
Thus edges start at rate $\alpha_i$ while the root is infected, and their delays are independent with density $q_{i+1}(v)/\chi_{i+1}$. Let $D$ be one such delay. If $u$ units of time have passed since it started and we are still waiting for its return, the probability that it returns in the next $h$ units of time is
\[
\Pp(D\le u+h\mid D>u) =\frac{\Pp(u<D\le u+h)}{\Pp(D>u)} =\frac{\int_u^{u+h}q_{i+1}(v)\,dv}{\int_u^\infty q_{i+1}(v)\,dv}.
\]

The root state has bounded, right-continuous paths, so $q_{i+1}$ is right-continuous. The instantaneous return rate is therefore
\[
\lim_{h\downarrow0}\frac{\Pp(D\le u+h\mid D>u)}{h} =\frac{q_{i+1}(u)}{\int_u^\infty q_{i+1}(v)\,dv} \le\frac1{\chi_{i+1}}.
\]

At the type-$i$ root, an edge $(s,w)$ is \emph{pending} at time $t$ if the root was infected at $s$ and $s\le t<w$. Let $A_j(t)$ be the event that the root is infected at $t$ with exactly $j$ edges pending. The expected time spent in this state is
\begin{equation}\label{eq:pending-occupation}
U_j=\int_0^\infty \Pp(A_j(t))\,dt.
\end{equation}

Count pending edges on the clock that runs only while the root is infected. A new edge raises the count by one. A return while the root is infected lowers it by one. If the root recovers with $j\ge1$ edges pending, the clock stops. Until the first of those edges returns, the root stays healthy, so no new edge enters the count. That return reinfects the root and leaves the other $j-1$ edges pending. The clock thus stops at count $j$ and resumes at $j-1$. A recovery while the root is infected with no pending edge is terminal, since the root then stays healthy forever. The expected number of edges started while the root is infected is $\alpha_i\chi_i<\infty$. Each infected period after the initial one begins at the return of a distinct such edge, so almost surely there is a last infected period. It ends at the next recovery mark. At that recovery no edge can be pending, since the first return would start another infected period. Thus exactly one recovery is terminal.

A recovery that finds the root infected with no edges pending is terminal, so
\[
1=\E[\text{number of terminal recoveries}] =\int_0^\infty \Pp(A_0(t))\,dt=U_0 .
\]

For $j\ge1$, on $A_{j-1}(t)$ an edge start gives an upward crossing at rate $\alpha_i$. On $A_j(t)$ a recovery gives a downward crossing at rate $1$, and each of the $j$ pending edges gives one at rate at most $1/\chi_{i+1}$. The count starts at zero and ends at zero at the terminal recovery, so each upward crossing from $j-1$ to $j$ is paired with a downward crossing. Thus
\begin{align*}
\alpha_i U_{j-1} &=\int_0^\infty \alpha_i\,\Pp(A_{j-1}(t))\,dt
   =\E[\text{number of crossings }j-1\to j]\\
 &=\E[\text{number of crossings }j\to j-1]\\
&\le \int_0^\infty \Pp(A_j(t))\,dt +\int_0^\infty \frac{j}{\chi_{i+1}}\,\Pp(A_j(t))\,dt =\left(1+\frac{j}{\chi_{i+1}}\right)U_j .
\end{align*}

For $j\ge1$, iterating the above bound gives
\[
U_j \ge \frac{\alpha_i}{1+j/\chi_{i+1}}\,U_{j-1} \ge\cdots\ge U_0\prod_{\ell=1}^j\frac{\alpha_i}{1+\ell/\chi_{i+1}} =\frac{\alpha_i^j}{\prod_{\ell=1}^j(1+\ell/\chi_{i+1})},
\]
since $U_0=1$. Each infected time belongs to exactly one level, so
\[
\chi_i=U_0+\sum_{j=1}^\infty U_j \ge 1+\sum_{j=1}^\infty \frac{\alpha_i^j}{\prod_{\ell=1}^j(1+\ell/\chi_{i+1})}. \qedhere
\]
\end{proof}

The terms for zero, one and two pending edges already suffice to show that $H_{\bb}(2/3)>1$.

\begin{proof}[Proof of the upper bound in Proposition~\ref{prop:upper-bounds}]
Take $a=2/3$ and suppose $H_{\bb}(a)\le1$, so all susceptibilities are finite. Write $t_i=\alpha_i/(1+1/\chi_{i+1}) =a^2b_i\chi_{i+1}^2/(\chi_{i+1}+1)$ for the $j=1$ term in
\eqref{eq:pending-levels} of Lemma~\ref{lem:pending-returns}.
Keeping the constant term and the first two summands gives
\begin{align*}
\chi_i &\ge 1+\frac{\alpha_i}{1+1/\chi_{i+1}}
       +\frac{\alpha_i^2}{(1+1/\chi_{i+1})(1+2/\chi_{i+1})}\\
&=1+t_i+\frac{1+1/\chi_{i+1}}{1+2/\chi_{i+1}}\,t_i^2 =1+t_i+\frac{\chi_{i+1}+1}{\chi_{i+1}+2}\,t_i^2.
\end{align*}

Since $t_i>0$, this gives $\chi_i>1$ for every $i$. In particular $(\chi_{i+1}+1)/(\chi_{i+1}+2)>2/3$, so
\[
\chi_i>1+t_i+\tfrac23t_i^2 =\frac{5+10t_i+(2t_i-1)^2}{6}\ge\frac{5+10t_i}{6}.
\]

Let $z_i=1/\chi_i$ and $\bar z=k^{-1}\sum_i z_i$. Multiplying $10t_i<6\chi_i-5$ by $z_i$ gives $10t_i z_i<6-5z_i$. By the definition of $t_i$, $a^2b_i=t_i z_{i+1}(1+z_{i+1})$. Since $\prod_i b_i=1$ and $z_{k+1}=z_1$,
\begin{align*}
(10a^2)^k &=\prod_{i=1}^k 10a^2b_i =\left(\prod_{i=1}^k 10t_i\right)
  \left(\prod_{i=1}^k z_{i+1}(1+z_{i+1})\right)\\
&=\prod_{i=1}^k (10t_i z_i)(1+z_i) <\prod_{i=1}^k (6-5z_i)(1+z_i).
\end{align*}

Taking $k$\textsuperscript{th} roots and applying the arithmetic--geometric mean inequality to the two products separately gives
\[
10a^2<\left(\prod_{i=1}^k(6-5z_i)\right)^{1/k} \left(\prod_{i=1}^k(1+z_i)\right)^{1/k} \le (6-5\bar z)(1+\bar z).
\]

Since $z_i=1/\chi_i$ and $H_{\bb}(a)=a^k\prod_i\chi_i\le1$,
\[
\bar z=\frac1k\sum_{i=1}^k z_i \ge \left(\prod_{i=1}^k z_i\right)^{1/k} =\frac{a}{H_{\bb}(a)^{1/k}}\ge a=\frac23 .
\]

The derivative of $(6-5z)(1+z)$ is $1-10z<0$ for $z\ge2/3$, so the preceding bounds give the contradiction
\[
10a^2<(6-5\bar z)(1+\bar z) \le \left(6-5\cdot\frac23\right)\left(1+\frac23\right) =\frac{40}{9}=10a^2 .
\]

Thus $H_{\bb}(2/3)>1=H_{\bb}(\kappa(\bb))$ by \eqref{eq:kappa}. Since $H_{\bb}$ is nondecreasing, $\kappa(\bb)<2/3$.
\end{proof}

\printbibliography

\subsection*{Statement on AI use}

The author used a large language model (ChatGPT Astra 6.2) during the preparation of this manuscript for editorial and expository support. The ideas of the proofs are due to the author and all AI-generated suggestions were reviewed, revised, and approved by the author, who takes full responsibility for the content of the manuscript.

\appendix
\section{Biregular trees with a large degree ratio}\label{app:biregular-ratio}

Let $m\le n$ be the offspring numbers of types $1$ and $2$, respectively. First let $m\to\infty$ with $n/m\to\theta$, so the normalised degree vector tends to $(\theta^{-1/2},\theta^{1/2})$. We ask how the limiting constant in Theorem~\ref{thm:periodic} behaves as $\theta\to\infty$.

\begin{proposition}\label{prop:biregular-ratio}
For $\theta\ge1$,
\[
\log\frac{(1+\sqrt\theta)^2}{1+2\sqrt\theta} \le \sqrt\theta\,\kappa(\theta^{-1/2},\theta^{1/2})^2 \le \log(1+\sqrt\theta).
\]

In particular, $\kappa(\theta^{-1/2},\theta^{1/2}) \sim\theta^{-1/4}\sqrt{(\log\theta)/2}$ as $\theta\to\infty$.
\end{proposition}

\begin{proof}
Take $a=\kappa(\theta^{-1/2},\theta^{1/2})$ and put $\gamma=\sqrt\theta\,a^2$. At this rate, Theorem~\ref{thm:periodic} gives
\[
a^2\chi_1\chi_2=1,\qquad \gamma\chi_1\chi_2=\sqrt\theta.
\]
Both susceptibilities are finite since each is at least $1$.

At a type-$1$ root, an edge starting at $s$ with delay $u$ has intensity $(a^2/\sqrt\theta)q_2(u)\,ds\,du$. Each infected period ends at a recovery, and each period after the first begins at the return of a distinct edge that started while the root was infected. Recoveries have rate one, so
\begin{align*}
\chi_1
 &=\E[\text{number of recoveries while the root is infected}]\\
 &\le 1+\E[\text{number of edges starting while the root is infected}]\\
&=1+\frac{a^2}{\sqrt\theta}\int_0^\infty q_1(s)\,ds \int_0^\infty q_2(u)\,du =1+\frac{a^2\chi_1\chi_2}{\sqrt\theta} =1+\frac1{\sqrt\theta}.
\end{align*}

At a type-$2$ root, edges $(s,s+u)$ have intensity $\gamma q_1(u)\,ds\,du$. Each edge represents a possible return through a type-$1$ branch. Starting the branch root infected, the probability it is infected at time $u$ is $q_1(u)$. Of this probability, $e^{-u}$ comes from having no recovery by $u$, and the rest comes from recovering and being infected again at $u$.

Independently keep each edge $(s,s+u)$ with probability $e^{-u}/q_1(u)$, and call the remaining edges extra edges. This splits the Poisson process into two independent Poisson processes with intensities $\gamma e^{-u}\,ds\,du$ for the kept edges, and $\gamma(q_1(u)-e^{-u})\,ds\,du$ for the excess.

Thus kept edges start at rate $\gamma$ and have independent delays with density $e^{-u}$ for $u>0$. Let $p(t)$ be the probability that the initially infected type-$2$ root is infected at $t$ using only kept edges. With the same recovery marks, every path using kept edges is also a path in the full process, so $p(t)\le q_2(t)$.

Use $A_j(t)$ and $U_j$ as in \eqref{eq:pending-occupation} for the process with only kept edges. Since $\gamma\int_0^\infty p(t)\,dt\le\gamma\chi_2<\infty$, the terminal recovery argument in the proof of Lemma~\ref{lem:pending-returns} gives $U_0=1$. Here a pending edge of age $u$ has return rate
\[
\frac{e^{-u}}{\int_u^\infty e^{-v}\,dv}=1.
\]

On $A_j(t)$ with $j\ge1$, the pending edges therefore return at total rate $j$. Recoveries contribute rate one, so the level-crossing calculation in that proof gives
\[
\gamma U_{j-1}=(j+1)U_j,\qquad U_j=\frac{\gamma^j}{(j+1)!}.
\]

Hence
\[
\int_0^\infty p(t)\,dt =\sum_{j=0}^\infty U_j =\frac1\gamma\sum_{j=0}^\infty\frac{\gamma^{j+1}}{(j+1)!} =\frac{e^\gamma-1}{\gamma}.
\]

The upper bound now follows from
\[
e^\gamma-1=\gamma\int_0^\infty p(t)\,dt \le\gamma\chi_2\le\gamma\chi_1\chi_2=\sqrt\theta.
\]

For the lower bound, any infection path from $0$ to $t$ either uses only kept edges or has a last extra edge $(s,w)$. To use $(s,w)$ as its last extra edge, the path must reach $s$ using either kind of edge and go from $w$ to $t$ using only kept edges. These two events depend on marks wholly within $[0,s]$ and $[w,t]$, respectively, so their probabilities multiply to $q_2(s)p(t-w)$. Revealing the extra edges in order of their endpoints as in Lemma~\ref{lem:effective-return}, an edge from $s$ to $w$ has intensity $\gamma(q_1(w-s)-e^{-(w-s)})\,ds\,dw$ independently of these marks. Counting all extra edges with the two required paths gives
\begin{align*}
q_2(t)
 &\le p(t)+\E[\text{number of such extra edges $(s,w)$}]\\
&=p(t)+\gamma\int_0^t\!\int_0^w q_2(s)\bigl(q_1(w-s)-e^{-(w-s)}\bigr)p(t-w)\,ds\,dw.
\end{align*}

With $u=w-s$ and $r=t-w$, Tonelli's theorem gives
\begin{align*}
\chi_2 &\le \int_0^\infty p(t)\,dt +\gamma\left(\int_0^\infty q_2(s)\,ds\right) \left(\int_0^\infty(q_1(u)-e^{-u})\,du\right)
         \left(\int_0^\infty p(r)\,dr\right)\\
&=\frac{e^\gamma-1}{\gamma} \bigl[1+\gamma\chi_2(\chi_1-1)\bigr].
\end{align*}

Dividing by $\chi_2$ and using $\gamma\chi_1\chi_2=\sqrt\theta$ gives
\begin{align*}
1 &\le(e^\gamma-1)\left(\frac1{\gamma\chi_2}+\chi_1-1\right)
   =(e^\gamma-1)\left(\frac{\chi_1}{\sqrt\theta}+\chi_1-1\right)\\
&\le(e^\gamma-1)\left(\frac{1+1/\sqrt\theta}{\sqrt\theta} +\frac1{\sqrt\theta}\right) =(e^\gamma-1)\frac{1+2\sqrt\theta}{\theta}.
\end{align*}

Together with the upper bound this gives
\[
\frac{(1+\sqrt\theta)^2}{1+2\sqrt\theta} =1+\frac{\theta}{1+2\sqrt\theta}\le e^\gamma\le 1+\sqrt\theta.
\]

The logarithms of the two sides differ by $\log((1+2\sqrt\theta)/(1+\sqrt\theta))<\log 2$. Hence,
\[
a^2=\frac{\gamma}{\sqrt\theta} =\frac{\log(1+\sqrt\theta)+O(1)}{\sqrt\theta} \sim\frac{\log\theta}{2\sqrt\theta}.\qedhere
\]
\end{proof}

\end{document}